\documentclass[12pt]{article}

\usepackage{amssymb,amsmath,graphicx}

\newcommand{\eps}{\varepsilon}

\newcommand{\lam}{\lambda}

\newcommand{\abs}[1]{\left|#1\right|}

\newcommand{\cB}{\mathcal{B}}

\newcommand{\bP}{\mathbf{Pr}}

\newtheorem{propn}{Proposition}
\newtheorem{cor}[propn]{Corollary}

\begin{document}

\begin{center}
{\large Improved bounds on the peak sidelobe level of binary sequences}
\end{center}

\begin{flushright}
Idris Mercer\\
Florida International University \\
\verb+imercer@fiu.edu+
\end{flushright}

\begin{abstract}
Schmidt proved in 2014 that if $\eps>0$,
almost all binary sequences of length~$n$
have peak sidelobe level between
$(\sqrt{2}-\eps)\sqrt{n\log n}$ and
$(\sqrt{2}+\eps)\sqrt{n\log n}$.
Because of the small gap between his upper
and lower bounds, it is difficult to find
improved upper bounds that hold for
almost all binary sequences. In this note,
we prove that if $\eps>0$,
then almost all binary sequences
of length~$n$ have peak sidelobe level
at most $\sqrt{2n(\log n-(1-\eps)\log\log n)}$,
and we provide a slightly better upper bound
that holds for a positive proportion of
binary sequences of length~$n$.
\end{abstract}

By a {\bf binary sequence} of {\bf length} $n$, we mean an $n$-tuple
$$
A = (a_0, a_1, \ldots, a_{n-1})
$$
where each $a_j$ is $+1$ or $-1$.
For $0 \le k \le n-1$, we define the (acyclic or aperiodic)
{\bf autocorrelations} of $A$ by
$$
c_k = \sum_{j=0}^{n-k-1} a_j a_{j+k}.
$$
Informally, $c_k$ measures how much the sequence $A$
resembles a version of itself that has been shifted by $k$ positions.

We let $\cB_n$ denote the set of all $2^n$ binary sequences of length~$n$.
For any $A \in \cB_n$, we have $c_0 = n$.
We refer to $c_1, \ldots, c_{n-1}$ as the {\bf nontrivial} autocorrelations of~$A$.
An old problem, arising in communications engineering but also of interest
as a stand-alone combinatorial problem, involves trying to find binary sequences
in $\cB_n$ whose nontrivial autocorrelations are `close' to zero in some sense.

For any $A \in \cB_n$, we define the {\bf peak sidelobe level} (PSL) of $A$ by
$$
\mu(A) = \max_{1\le k\le n-1} \abs{c_k}.
$$
We consider $A$ to be a `good' sequence if $\mu(A)$ is small.
If $A$ is a constant sequence, then trivially $\mu(A)=n-1$,
but very informally speaking, if $A$ is `random' then $\mu(A)$
tends to be significantly smaller than~$O(n)$.
Many authors have investigated upper bounds for $\mu(A)$.
(For an excellent survey, see~\cite{JY}.)
We might try to find upper bounds for $\mu(A)$ that hold
for some sequences $A\in\cB_n$,
or that hold for almost all sequences $A\in\cB_n$.

To make this more precise, we turn $\cB_n$ into a probability space
by supposing the $a_j$ are independent Rademacher variables
(i.e., random variables each equally likely to be $+1$ or $-1$).
This is equivalent to assigning equal weight to each of the
$2^n$ sequences in~$\cB_n$, and for any function $f(n)$,
the probability that $\mu(A) \le f(n)$ is equal to the proportion
of sequences $A\in\cB_n$ that satisfy $\mu(A) \le f(n)$.
We say $\mu(A) \le f(n)$ for `almost all' sequences $A\in\cB_n$
if
$$
\lim_{n\to\infty} \bP\big[\mu(A) \le f(n)\big] = 1.
$$
We also define
$$
\mu_{\min}(n) = \min_{A\in\cB_n} \mu(A)
$$
so then if $\mu(A) \le f(n)$ for a nonzero proportion of
sequences $A\in\cB_n$, we have $\mu_{\min}(n) \le f(n)$.

In 2014, Schmidt proved \cite{Sch14} (improving upon previous results by
Alon, Litsyn \& Shpunt \cite{ALS}, the current author \cite{Mer},
and Moon \& Moser \cite{MM}) that if we fix $\eps>0$,
then the probability
\begin{equation} \label{Schmidt}
\bP\Big[(\sqrt{2}-\eps)\sqrt{n\log n} \le \mu(A) \le (\sqrt{2}+\eps)\sqrt{n\log n}\Big]
\end{equation}
approaches 1 as $n$ approaches infinity
(informally, almost all sequences $A\in\cB_n$
have peak sidelobe level `close' to $\sqrt{2n\log n}$).
Here and throughout this article, `log' means natural log.

Earlier, Schmidt \cite{Sch12} gave an explicit construction
showing that for each $n>1$, there is a sequence $A\in\cB_n$
satisfying $\mu(A)\le\sqrt{2n\log(2n)}$.
He also gave numerical evidence for the conjecture that
his sequences satisfy $\mu(A)=O(\sqrt{n\log\log n})$.
As pointed out in~\cite{JY}, several authors have conjectured
that there is an infinite family of binary sequences
satisfying $\mu(A)=O(\sqrt{n})$, but this has not been proved.
In fact, the best upper bounds that have been proved to hold
{\it either} for a positive proportion of sequences
{\it or} for almost all sequences appear to be
of the form $\mu(A)=O(\sqrt{n\log n})$.

Because of the lower bound in (\ref{Schmidt}),
it is not possible to prove that almost all sequences $A\in\cB_n$
satisfy an upper bound of the form $\mu(A)=o(\sqrt{n\log n})$.
However, if $f(n)$ is a certain function of~$n$
that approaches infinity more slowly than $\log n$,
we can prove that almost all sequences $A\in\cB_n$
satisfy $\mu(A)\le\sqrt{2n(\log n-f(n))}$.
By slightly modifying $f(n)$, we can find a similar upper bound
that holds for a positive proportion of sequences in~$\cB_n$.

More specifically, we prove the following proposition and corollaries.

\begin{propn}
Let $\psi(n)$ be a function of $n$.
(The conclusion is interesting only if $\psi(n)$
approaches infinity with~$n$.)
Then the proportion of sequences $A \in \cB_n$
satisfying
$$
\mu(A) > \sqrt{2n\psi(n)}
$$
is bounded above by
$$
\frac{2n}{\psi(n)e^{\psi(n)}}.
$$
\end{propn}

\begin{cor}
Let $\eps>0$.
Then the proportion of sequences $A \in \cB_n$ satisfying
$$
\mu(A) > \sqrt{2n\big(\log n - (1-\eps)\log\log n\big)}
$$
approaches 0 when $n$ approaches infinity.
\end{cor}

\begin{cor}
Let $\eps>0$.
Then the proportion of sequences $A \in \cB_n$ satisfying
$$
\mu(A) > \sqrt{2n\big(\log n - \log\log n + \log 2 + \eps \big)}
$$
is strictly less than 1 for all sufficiently large~$n$.
\end{cor}

Notice that Corollary 2 says that
$$
\mu(A) \le \sqrt{2n\big(\log n - (1-\eps)\log\log n\big)}
$$
for almost all sequences $A\in\cB_n$, and Corollary 3
says that eventually,
$$
\mu_{\min}(n) \le \sqrt{2n\big(\log n - \log\log n + \log 2 + \eps \big)}.
$$
Note that $\log2 \approx 0.693$. In addition to having a bound
for $\mu_{\min}(n)$ that holds for all sufficiently large~$n$,
it may be of interest to have a bound that holds for all $n>1$.

\begin{cor}
For all $n>1$, the proportion of sequences $A\in\cB_n$ satisfying
$$
\mu(A) > \sqrt{2n\big(\log n - \log\log n + 0.862\big)}
$$
is strictly less than 1.
\end{cor}

Notice that Corollary 4 says that
$$
\mu_{\min}(n) \le \sqrt{2n\big(\log n - \log\log n + 0.862\big)}
$$
for all $n>1$.

{\bf Proof of Proposition 1:}

As mentioned before,
we turn $\cB_n$ into a probability space by supposing the $a_j$ to be
independent Rademacher variables,
which is equivalent to assigning equal weight to all $2^n$ sequences in~$\cB_n$.

Note that the autocorrelation
$$
c_k = a_0 a_k + a_1 a_{k+1} + \cdots + a_{n-k-1} a_{n-1}
$$
is a sum of $n-k$ terms, each of which is $\pm1$.
In fact, those $n-k$ terms are independent. (This is straightforward
but not quite trivial; for a proof, see~\cite{Mer}.)
If $1 \le k \le n-1$, then $c_{n-k}$ is a sum of $k$ independent
Rademacher variables, so we can use Chernoff-type bounds
(see, e.g., Corollary A.1.2 in Appendix A of \cite{AS})
to conclude that if $\lam>0$, then
$$
\bP\big[\abs{c_{n-k}}>\lam\big] < 2\exp(-\lam^2/2k).
$$
If $\lam = \sqrt{2n\psi(n)}$, this becomes
$$
\bP\Big[\abs{c_{n-k}}>\sqrt{2n\psi(n)}\Big] < 2\exp(-n\psi(n)/k).
$$
We call a sequence $A \in \cB_n$ `good' if $\mu(A) \le \sqrt{2n\psi(n)}$,
and `bad' otherwise.
Then $A$ is bad if and only if $\abs{c_{n-k}} > \sqrt{2n\psi(n)}$
for some $k=1,\ldots,n-1$.
An overestimate for $\bP[\mbox{$A$ is bad}]$ is
$$
\sum_{k=1}^{n-1} \bP\Big[\abs{c_{n-k}}>\sqrt{2n\psi(n)}\Big] < \sum_{k=1}^{n-1} 2\exp(-n\psi(n)/k).
$$
Now, consider the function
$$
g(x) = 2\exp(-\psi(n)/x)
$$
on the interval $x \in [\frac1n, 1]$.
The function $g(x)$ is an increasing function of $x$ on that interval,
so a left-endpoint Riemann sum will be an underestimate for an integral:
\begin{align*}
\sum_{k=1}^{n-1} g\Big(\frac{k}{n}\Big) \frac1n & < \int_{1/n}^1 g(x)dx \\[1ex]
\Longrightarrow \sum_{k=1}^{n-1} g\Big(\frac{k}{n}\Big) & < n\int_{1/n}^1 g(x)dx \\[1ex]
\Longrightarrow \sum_{k=1}^{n-1} 2\exp(-n\psi(n)/k) & < 2n \int_{1/n}^1 \exp(-\psi(n)/x)dx \\[1ex]
\Longrightarrow \bP[\mbox{$A$ is bad}] & < 2n \int_{1/n}^1 \exp(-\psi(n)/x)dx.
\end{align*}
We will now perform the substitution $u=\psi(n)/x$ on this integral. We have
\begin{align*}
u & = \psi(n)x^{-1} \\
du & = -\psi(n)x^{-2}dx \\ 
-\big(x^2/\psi(n)\big)du & = dx \\
x=1/n & \Rightarrow u = n\psi(n) \\
x=1 & \Rightarrow u = \psi(n) \\
x & = \psi(n)/u \\
x^2 & = \big(\psi(n)\big)^2/u^2 \\
x^2 / \psi(n) & = \psi(n)/u^2 \\
dx = -\big(x^2/\psi(n)\big)du & = -\big(\psi(n)/u^2\big)du
\end{align*}
and so the above integral becomes
\begin{align*}
2n \int_{1/n}^1 \exp(-\psi(n)/x)dx
& = 2n \int_{n\psi(n)}^{\psi(n)} \exp(-u) \Big(-\frac{\psi(n)}{u^2}\Big)du \\[1ex]
& = 2n\psi(n) \int_{\psi(n)}^{n\psi(n)} \frac{1}{u^2e^u}du.
\end{align*}
That is, we have
$$
\bP[\mbox{$A$ is bad}] < 2n\psi(n) \int_{\psi(n)}^{n\psi(n)} \frac{1}{u^2e^u}du.
$$
Now since the function $h(u) = 1/u^2e^u$ decreases very rapidly,
a rather crude upper bound will suffice. We have
$$
\int_{\psi(n)}^{n\psi(n)} \frac{1}{u^2e^u}du < \int_{\psi(n)}^\infty \frac{1}{u^2e^u}du.
$$
On the interval $u\in[\psi(n),\infty)$, we have $u^2 > (\psi(n))^2$, so we have
$$
\int_{\psi(n)}^\infty \frac{1}{u^2e^u}du < \frac{1}{(\psi(n))^2} \int_{\psi(n)}^\infty e^{-u} du
= \frac{1}{(\psi(n))^2} e^{-\psi(n)}.
$$
This implies that we have
$$
\bP[\mbox{$A$ is bad}] < 2n\psi(n) \cdot \frac{1}{(\psi(n))^2} e^{-\psi(n)}
= \frac{2n}{\psi(n)e^{\psi(n)}},
$$
completing the proof of Proposition 1.

{\bf Proof of Corollary 2:}

Let $\eps>0$, and define
$$
\psi(n) = \log n - (1-\eps) \log\log n.
$$
By Proposition 1, the proportion of sequences $A\in\cB_n$ satisfying
$\mu(A) > \sqrt{2n\psi(n)}$ is bounded above by
$$
\frac{2n}{\psi(n)e^{\psi(n)}}.
$$
Observe that for this choice of $\psi(n)$, we have
\begin{align*}
\exp\big(\psi(n)\big) & = \exp(\log n) \exp\big(-(1-\eps)\log\log n\big) \\
& = n \exp\Big(\log\big((\log n)^{-(1-\eps)}\big)\Big) \\
& = n(\log n)^{-(1-\eps)}
\end{align*}
which means that we have
\begin{align*}
\frac{2n}{\psi(n)e^{\psi(n)}}
& = \frac{2n}{\psi(n)n(\log n)^{-(1-\eps)}}
= \frac{2(\log n)^{1-\eps}}{\psi(n)} \\[1ex]
& = \frac{2(\log n)^{1-\eps}}{\log n - (1-\eps)\log\log n},
\end{align*}
which approaches 0 as $n$ approaches infinity.

{\bf Proof of Corollary 3:}

Let $\eps>0$, and define
$$
\psi(n) = \log n - \log\log n + \log 2 + \eps.
$$
By Proposition 1, the proportion of sequences $A\in\cB_n$ satisfying
$\mu(A) > \sqrt{2n\psi(n)}$ is bounded above by
$$
\frac{2n}{\psi(n)e^{\psi(n)}}.
$$
Observe that for this choice of $\psi(n)$, we have
\begin{align*}
\exp\big(\psi(n)\big) & = \exp(\log n) \exp(-\log\log n) \exp(\log2) \exp(\eps) \\
& = n (\log n)^{-1} 2\exp(\eps)
\end{align*}
which means that we have
\begin{align*}
\frac{2n}{\psi(n)e^{\psi(n)}}
& = \frac{2n}{\psi(n)n(\log n)^{-1}2\exp(\eps)}
= \frac{\log n}{\exp(\eps)\psi(n)} \\[1ex]
& = \frac{\log n}{\exp(\eps)(\log n - \log\log n + \log2 + \eps)},
\end{align*}
which approaches $1/\exp(\eps)<1$ as $n$ approaches infinity.

To prove Corollary 4, we use the following fact.

{\bf Fact.} If $n>1$ and $K$ is a constant, then
$$
\frac{K-\log\log n}{\log n} \ge \frac{-1}{e^{K+1}}.
$$
{\bf Proof of Fact.} Consider the function
$$
f(x) = \frac{K-\log x}{x}
$$
for $x>0$. Using elementary calculus, we find
$$
f'(x) = \frac{\log x -(K+1)}{x^2}
$$
which is negative when $0<x<e^{K+1}$ and positive when $x>e^{K+1}$.
It follows that for all $x>0$, we have
$$
f(x) \ge f(e^{K+1}) = \frac{-1}{e^{K+1}}
$$
and therefore for all $n>1$, we have
$$
\frac{K-\log\log n}{\log n} = f(\log n) \ge \frac{-1}{e^{K+1}}.
$$

{\bf Proof of Corollary 4:}

Suppose $n>1$, and define
$$
\psi(n) = \log n - \log\log n + 0.862.
$$
By Proposition 1, the proportion of sequences $A\in\cB_n$
satisfying $\mu(A) > \sqrt{2n\psi(n)}$ is bounded above by
$$
\frac{2n}{\psi(n)e^{\psi(n)}}.
$$
Observe that for this choice of $\psi(n)$, we have
\begin{align*}
\exp\big(\psi(n)\big) & = \exp(\log n) \exp(-\log\log n) \exp(0.862) \\
& = n (\log n)^{-1} e^K
\end{align*}
where for brevity, we write $K=0.862$. We then have
\begin{align*}
\psi(n)e^{\psi(n)} & = \big( \log n - \log\log n + K \big) \cdot n(\log n)^{-1}e^K \\
& = e^K \Big( 1 + \frac{K-\log\log n}{\log n} \Big) n
\end{align*}
and then the fact stated earlier implies
$$
\psi(n)e^{\psi(n)} \ge e^K \Big( 1 + \frac{-1}{e^{K+1}} \Big) n = \Big( e^K - \frac1e \Big) n.
$$
Now note that
$$
e^K - \frac1e = e^{0.862} - \frac1e > 2.00001
$$
so we have
$$
\frac{2n}{\psi(n)e^{\psi(n)}} < \frac{2n}{2.00001n} = \frac{2}{2.00001} < 1
$$
which completes the proof of the corollary.

Finally, to illustrate how the bound $\mu(A) \le \sqrt{2n\big(\log n - \log\log n + 0.862\big)}$
compares to the bound $\mu(A) \le \sqrt{2n\log n}$, we list numerical values of these bounds
for several large values of $n$.

\begin{center}
\begin{tabular}{c|c|c}
$n$ & $\sqrt{2n\log n}$ & $\sqrt{2n(\log n-\log\log n+0.862)}$ \\ \hline
1000 & 117.54 & 108.85 \\
2000 & 174.37 & 160.43 \\
3000 & 219.18 & 201.81 \\
4000 & 257.59 & 237.33 \\
5000 & 291.84 & 269.02 \\
6000 & 323.10 & 297.96 \\
7000 & 352.07 & 324.79 \\
8000 & 379.20 & 349.93 \\
9000 & 404.83 & 373.69 \\
10000 & 429.19 & 396.28
\end{tabular}
\end{center}

\bibliographystyle{amsplain}

\end{document}